\documentstyle[12pt,amsfonts,amssymb,leqno]{article}

\newfont\got{eufm10}

\newtheorem{proposition}{Proposition}[section]
\newtheorem{thm}[proposition]{Theorem}

\newtheorem{lemma}[proposition]{Lemma}

\newcounter{secnum}

\renewcommand{\thefootnote}{\fnsymbol{footnote}}
\renewcommand{\thefootnote}{arabic{footnote}}
\renewcommand{\thefootnote}{\fnsymbol{footnote}}

\newcommand{\zerohandgrenade}{0^{\: \mid^{\! \! \! \bullet}}}

\begin{document}
\setcounter{section}{+1}

\begin{center}
{\Large \bf More on mutual stationarity}
\end{center}
\begin{center}
{\fbox{Preliminary version}}
\end{center}
\begin{center}
{\large Ralf Schindler}
\renewcommand{\thefootnote}{arabic{footnote}}
\end{center}
\begin{center} 
{\footnotesize
{\it Institut f\"ur Formale Logik, Universit\"at Wien, 1090 Wien, Austria}} 
\end{center}

\begin{center}
{\tt rds@logic.univie.ac.at}

{\tt http://www.logic.univie.ac.at/${}^\sim$rds/}\\
\end{center}

\begin{abstract}
\noindent
We prove that in the core model below $\zerohandgrenade$ the following holds true.

\noindent Let $k \in OR \setminus 1$. There is then a sequence 
$(S^n_i \ \colon \ i \in OR \setminus k , 0 < n < \omega)$
such that 
 
\noindent $\bullet \ $ for all $i \in OR \setminus k$ 
and $n < \omega \setminus 1$ do we have that 
$S^n_i \subset \aleph_{i+1}$ is 
stationary in $\aleph_{i+1}$ and $\alpha 
\in S^n_i \Rightarrow cf(\alpha) = \aleph_k$, and

\noindent $\bullet \ $ for all limit ordinals $\lambda$ and for 
all $f \colon \lambda \rightarrow 
\omega \setminus 1$ 
do we have that $(S^{f(i)}_i \ \colon \ k \leq i < \lambda)$ is mutually stationary
if and only if the range of $f$ is finite.
\end{abstract}

Let $A$ be a non-empty set of regular cardinals. Recall that $(S_\kappa \colon
\kappa \in A)$ is called mutually stationary if and only if for each large enough
regular cardinal $\theta$ and for each model ${\frak A}$
with universe $H_\theta$ there is some $X \prec {\frak A}$ such that 
for all $\kappa \in
X \cap A$, $sup(X \cap \kappa) \in S_\kappa$. We refer the reader to \cite{fm} (in
particular, to \cite[Section 7]{fm}). It is an  
open problem to decide whether there is a model of set theory in which
$(S_n \colon n<\omega)$ must be mutually stationary provided 
each individual $S_n \subset
\aleph_n$ is stationary.

Foreman and Magidor have shown that in $L$ there is some $(S_n \colon n<\omega)$ such
that each $S_n \subset
\aleph_n$ is stationary, but $(S_n \colon n<\omega)$ is not mutually stationary
(cf. \cite[Theorems 24 and 27]{fm}). The purpose of this paper is to extend their
result \cite[Theorem 27]{fm} to the core model below 
$\zerohandgrenade$ which we introduced in
\cite{habil}.

We shall need the following lemma, which is essentially due to Baumgartner
\cite{baum}.   

\begin{lemma}\label{variations}
Let $m < \omega$, and
let $A_0$, $A_1$, ..., $A_m$ be non-empty sets of successor cardinals such that 
$sup(A_l) < min(A_{l+1})$ for all $l<m$. 
Let $(S_\kappa \colon \kappa \in A_l)$ 
be mutually stationary for all $l \leq m$ and such that
$\forall \kappa \in A_l \ \forall \alpha \in S_\kappa \ cf(\alpha) <
min(A_0)$. 
Then $(S_\kappa \colon \kappa \in \bigcup_{l \leq m} A_l)$ is 
mutually stationary as well.
\end{lemma}

{\sc Proof.} It certainly suffices to prove \ref{variations} for $m=1$.
Let $\nu$ be the cardinal predecessor of $min(A_0)$.
Let ${\frak A}$ be a model expanding $(H_\theta;\in)$ for some large regular
$\theta$. As $(S_\kappa \colon \kappa \in A_1)$ 
is mutually stationary we may pick some $X \prec {\frak A}$ such that $sup(A_0)
\subset X$ and for all $\kappa \in X \cap A_1$, $sup(\kappa \cap X) \in S_\kappa$.
Pick $F \colon (X \cap A_1) \times \nu \rightarrow X$ such that for all
$\kappa \in X \cap A_1$, $F(\kappa,-)$ is cofinal in $sup(\kappa \cap X)$.
Expand $X$ by $F$ to get $(X,F)$. As $(S_\kappa \colon \kappa \in A_0)$ 
is mutually stationary we may pick some $Y \prec (X,F)$ such that $\nu \subset Y$ and
for all $\kappa \in Y \cap A_0$, $sup(\kappa \cap Y) \in S_\kappa$. 
However, due to the presence of $F$, we'll also have that
for all $\kappa \in Y \cap A_1$, $sup(\kappa \cap Y) = sup(\kappa \cap X) \in
S_\kappa$. As $Y \prec {\frak A}$, we are done.
$\square$ {\scriptsize Lemma \ref{variations}}

\bigskip
We can now state and prove our main result. Our proof will 
closely follow \cite[Section 7.2]{fm} to a certain extent.

\begin{thm}\label{main-thm}
Suppose that $V=K$, where $K$ is the core model below $\zerohandgrenade$
from \cite{habil}. 
Let $k \in OR \setminus 1$.
There is then a sequence 
$(S^n_i \ \colon \ i \in OR \setminus k , 0 < n < \omega)$
such that 
 
\noindent $\bullet \ $ for all $i \in OR \setminus k$ 
and $n < \omega \setminus 1$ do we have that 
$S^n_i \subset \aleph_{i+1}$ is 
stationary in $\aleph_{i+1}$ and $\alpha 
\in S^n_i \Rightarrow cf(\alpha) = \aleph_k$, and

\noindent $\bullet \ $ for all limit ordinals $\lambda$ and for 
all $f \colon \lambda \rightarrow 
\omega \setminus 1$ 
do we have that $(S^{f(i)}_i \ \colon \ k \leq i < \lambda)$ is mutually stationary
if and only if the range of $f$ is finite.
\end{thm}

{\sc Proof.} 
We commence by defining the sequence 
$(S^n_i \ \colon \ i \in OR \setminus k, 0 < n < \omega)$.
If $\alpha$ is a singular ordinal, then we let 
$(\gamma(\alpha),n(\alpha),\delta(\alpha),{\vec p}(\alpha))$ be the lexicographically
least tuple $(\gamma,n,\delta,{\vec p})$ such that 
$$Hull_n^{K||\gamma}(\delta \cup {\vec p}) \cap \alpha$$ is cofinal in $\alpha$.
For an ordinal $i \geq k$ and a natural number $n$ we define
$$S^n_i = \{ \alpha < \aleph_{i+1} \ \colon \ cf(\alpha) = \aleph_k \wedge n(\alpha) =
n \}.$$

We are now going to show that $(S^n_i \ \colon \ i \in OR \setminus k, 0 < 
n < \omega)$ witnesses
the truth of Theorem \ref{main-thm}.

We shall first prove that if $\lambda$ is a limit ordinal and $f \colon \lambda
\rightarrow \omega \setminus 1$ is such that $ran(f)$ is finite then
$(S^{f(i)}_i \ \colon \ k \leq
i < \lambda)$ is mutually stationary. By Lemma \ref{variations}
it will be enough if we prove this under the assumption that $ran(f) = \{ n \}$ for
some $n \in \omega \setminus 1$.

Let ${\frak A}$ be a model with universe $K|\kappa$,
where $\kappa$ is a large regular cardinal. Let $\gamma_0$ be least such that
${\frak A} \in K|\gamma_0$, and let $\gamma$ be the $\aleph_k^{th}$ ordinal $\beta$
such that $\beta > \gamma_0$ and $K|\beta \prec_{\Sigma_n} K|\kappa^+$. Let $$X =
Hull_{n}^{K|\gamma}(\aleph_k \cup \{ {\frak A} \}).$$ One can then argue as for 
\cite[Lemma 25]{fm} that for every $i \in X \cap \lambda$, $X \cap \aleph_{i+1} \in
S^{n}_i$. The only thing to notice that the use of the condensation lemma for $L$
can be replaced by \cite[\S 8, Lemma 4]{ronald} in a straightforward way.

We shall now prove that if $\lambda$ is a limit ordinal and $f \colon \lambda
\rightarrow \omega \setminus 1$ is such that
$(S^{f(i)}_i \ \colon \ k \leq i < \lambda)$ is mutually stationary then
$ran(f)$ is finite.

Let $\kappa$ be a large regular cardinal,
and let $$N \prec (K||\kappa;\in,...)$$ be such that
for all $i \in N \cap \lambda$ do 
we have that $sup(N \cap \aleph_{i+1}) \in S^{f(i)}_i$.
In particular $cf(sup(N \cap \aleph_{i+1})) = \aleph_k > \omega$ for each
such $i$. Let $$\pi \ \colon \ {\bar K} \cong N \prec K||\kappa$$ be such that
${\bar K}$ is transitive, and let ${\bar \lambda} = \pi^{-1}{\rm " }\lambda$. 
Let $\delta = c.p.(\pi)$, and let $\eta \leq \delta$ be least such that there is
some $X \in ({\cal P}(\eta) \cap K) \setminus {\bar K}$. Note that 
$({\cal P}(\delta) \cap K) \setminus {\bar K} \not= \emptyset$, because otherwise 
$\pi \upharpoonright {\cal P}(\delta) \cap K$ would be an extender which collapses
the cardinal $\pi(\delta)$.

Let $({\cal U},{\cal T})$ denote the (padded) coiteration
of ${\bar K}$ with $K$. We'll have that $[0,\infty)_U \cap {\cal D}^{\cal U} =
\emptyset$, and that $\pi^{\cal U}_{0 \infty} \upharpoonright \delta = id$. 

We first want to show that ${\cal U}$ is trivial (i.e., that
${\bar K}$ doesn't move in the comparison with $K$), and that if $\vartheta$ is least
such that ${\cal M}_\vartheta^{\cal T} \trianglerighteq {\bar K}||{\bar \lambda}$ then
${\cal M}_\vartheta^{\cal T}$ is set-sized and
$\rho_\omega({\cal M}_\vartheta^{\cal T}) < {\bar \lambda}$.

\bigskip
{\bf Claim 1.} Let $\mu$ be such that $\pi^{\cal U}_{0 \infty} \upharpoonright \mu =
id$. Then for no $F = E^{\cal T}_\nu$ do we have that ${\bar \mu} = c.p.(F) < \eta$ and
$F({\bar \mu}) \leq \mu$.

\bigskip
{\sc Proof.} Otherwise $\pi \circ F$ would be an extender which collapses
the cardinal $\pi \circ F(c.p.(F))$.\footnote{This argument heavily uses that our
premice are indexed \`a la Jensen rather than \`a la Mitchell-Steel.} 
$\square$ {\scriptsize Claim 1} 

\bigskip
{\bf Claim 2.} ${\cal U}$ is trivial.

\bigskip
{\sc Proof.} Assume not. Let $F = E_\epsilon^{\cal U}$,
where $\epsilon+1$ is least in $(0,\infty]_U$,
and let $\mu = c.p.(F) \geq \delta$. Let $\beta 
< lh({\cal U}) = lh({\cal T})$
be minimal with ${\cal M}_\beta^{\cal T} \trianglerighteq {\bar K}|\mu^{+{\bar K}}$. 
We let $(\kappa_\gamma \colon
\gamma \leq \theta)$ enumerate the cardinals of ${\bar K}$ in the
half-open interval
$[\eta,\mu^{+{\bar K}})$, and we let $\lambda_\gamma = \kappa_\gamma^{+{\bar K}}$ for
$\gamma \leq \theta$.
For each $\gamma \leq \theta$
we let $\delta(\gamma) \leq \beta$ be the least $\delta$ such that 
${\cal M}^{\cal T}_{\delta} | \lambda_\gamma = {\bar K} | \lambda_\gamma$, and
we let ${\cal P}_\gamma$ be the largest initial segment ${\cal P}$ of 
${\cal M}^{\cal T}_{\delta(\gamma)}$ 
such that all bounded subsets of $\lambda_\gamma$ which are in ${\cal P}$ are in 
${\bar K}$ as well.
By Claim 1 we shall have that $\lambda_\gamma = \kappa_\gamma^{+{\cal P}_\gamma}$ and
$\rho_\omega({\cal P}_\gamma) \leq \kappa_\gamma$ for all $\gamma \leq \theta$.
Moreover, ${\cal P}_\gamma$ is sound above $\kappa_\gamma$.

We let ${\vec {\cal P}}$ denote the
phalanx $$((K^\frown({\cal P}_\gamma \colon \gamma \leq \theta)^\frown {\cal M}^{\cal
T}_\epsilon),\eta^\frown(\lambda_\gamma \colon
\gamma \leq \theta)).$$
Let ${\cal U}^*$ and ${\cal V}$ be the padded iteration trees
arising from the comparison of ${\vec {\cal P}}$ with $K$. 
We understand that ${\cal U}^*$ either has a last ill-founded model, or else
that ${\cal M}_\infty^{{\cal U}^*}$ and ${\cal M}_\infty^{{\cal V}}$ are lined up.
The following says that it is the latter which will hold. 

\bigskip
{\bf Subclaim 1.} ${\vec {\cal P}}$ is coiterable with $K$.

\bigskip
{\sc Proof.} 
Suppose that ${\cal U}^*$ has a last ill-founded model. Let $\sigma \colon {\bar
H} \rightarrow H_\Omega$ be such that $\Omega$ is regular and large enough,
${\bar H}$ is countable and transitive, and $\{ {\vec {\cal P}} , {\cal U}^* \}
\subset ran(\sigma)$. Then $\sigma^{-1}({\cal U}^*)$ witnesses that 
$\sigma^{-1}({\vec {\cal P}})$ is not iterable. 

Fix $\gamma \in \theta+1 \cap
ran(\sigma)$ for a moment.
Let
$${\cal Q}_\gamma = Ult_n(\sigma^{-1}({\cal P}_\gamma);\sigma \upharpoonright 
\sigma^{-1}({\cal P}_\gamma | \lambda_\gamma)) {\rm , }$$
where $n<\omega$ is such that $\rho_{n+1}({\cal P}_\gamma) \leq \kappa_\gamma
< \rho_{n}({\cal P}_\gamma)$. Let $\sigma_\gamma \supset \sigma \upharpoonright 
\sigma^{-1}({\cal P}_\gamma | \lambda_\gamma)$ denote the canonical
embedding from $\sigma^{-1}({\cal P}_\gamma)$ into ${\cal Q}_\gamma$.
Notice that by $cf(\lambda_\gamma) > \omega$ it is clear that $sup \ \sigma{\rm "} 
\sigma^{-1}(\lambda_\gamma)$ is not cofinal in $\lambda_\gamma$.
Hence if $k_\gamma$ denotes the canonical embedding from ${\cal Q}_\gamma$ into
${\cal P}_\gamma$ then $sup \ \sigma{\rm "} 
\sigma^{-1}(\lambda_\gamma) = \kappa_\gamma^{+{\cal Q}_\gamma}
= k_\gamma^{-1}(\lambda_\gamma)$ is the critical point
of $k_\gamma$. 

Unfortunately, 
${\cal Q}_\gamma$ might not be premouse but rather a proto-mouse; 
this will in fact be
the case if ${\cal P}_\gamma$ has a top extender with critical point $\mu \leq
\kappa_\gamma$ and $n=1$, 
as then $\sigma_\gamma$ is discontinuous at 
$\sigma_\gamma^{-1}(\mu^{+{\cal P}_\gamma})$. 
Let us therefore define an
object ${\cal R}_\gamma$ as follows. 
We set ${\cal R}_\gamma = {\cal Q}_\gamma$ if
${\cal Q}_\gamma$ is a premouse. Otherwise, if $\mu \leq \kappa_\gamma$ 
is the critical
point of the top extender of ${\cal P}_\gamma$,
we let $${\cal R}_\gamma = Ult_{n}({\cal Q}_\gamma || \rho;G){\rm , }$$
where $G$ is the top extender fragment of ${\cal Q}_\gamma$, $\rho$ is maximal with
$$\mu^{+{{\cal Q}_\gamma || \rho}} = sup \ \sigma_\gamma{\rm " }
\sigma^{-1}(\mu)^{+{\sigma^{-1}({\cal P}_\gamma)}}{\rm , }$$
and $n<\omega$ is such that $\rho_{n+1}({\cal Q}_\gamma) \leq \mu < 
\rho_{n}({\cal Q}_\gamma)$.

Now because ${\cal R}_\gamma$ and ${\cal P}_\gamma$
are both sound above $\kappa_\gamma$ 
we can apply \cite[\S 8, Lemma 4]{ronald} (or an elaboration of it in case that 
${\cal Q}_\gamma$ is a proto-mouse) and deduce that
${\cal R}_\gamma \trianglelefteq {\cal P}_\gamma$.
But as $\kappa_\gamma^{+{\cal R}_\gamma} < \lambda_\gamma = 
\kappa_\gamma^{+{\cal P}_\gamma}$, this means that ${\cal R}_\gamma \triangleleft
{\cal P}_\gamma | \lambda_\gamma = {\cal M}^{\cal U}_\epsilon | \lambda_\gamma$. 
Therefore ${\cal
R}_\gamma$ is an initial segment of ${\cal M}^{\cal U}_\epsilon$.
 
Now if ${\vec {\cal R}}$ denotes the phalanx
$$(K^\frown({\cal R}_\gamma \colon \gamma \in \theta+1 \cap
ran(\sigma))^\frown{\cal M}_\epsilon^{\cal T},\eta^\frown(\lambda_\gamma \colon 
\gamma \in \theta+1 \cap
ran(\sigma)){\rm , }$$
then,
due to the existence of the family of maps
$$\sigma \upharpoonright \sigma^{-1}(K), (\sigma_\gamma \colon \gamma \in \theta+1 \cap
ran(\sigma)), \sigma \upharpoonright 
\sigma^{-1}({\cal M}_\epsilon^{\cal T})$$
we know that ${\vec {\cal R}}$ cannot be iterable, as $\sigma^{-1}({\vec {\cal P}})$
is not iterable. However, any iteration 
of ${\vec {\cal R}}$ can be construed as an iteration of
$((K,{\cal M}_\epsilon^{\cal U}),\delta)$, and thus in turn of
$((K,{\bar K}),\delta)$. But $((K,{\bar K}),\delta)$ is iterable. Contradiction!
$\square$ {\scriptsize Subclaim 1} 

\bigskip
Now notice that ${\cal U}^*
\upharpoonright \epsilon$ is trivial, ${\cal V} \upharpoonright \epsilon =
{\cal T} \upharpoonright \epsilon$, and that 
$F = E^{\cal U}_\epsilon = E^{{\cal U}^*}_\epsilon$ will be the first
extender used in ${\cal U}^*$. By \cite[Lemma 2.7]{habil}, 
we'll then\footnote{This is the only place in this paper 
where we really use the assumption that $K$ is below $\zerohandgrenade$.
If the assumption that $K$ is below $\zerohandgrenade$ is dropped
then at the time of writing I don't see how to prove the iterability of the phalanx
needed to verify Claim 2.} in fact 
have that no extender from ${\cal U}^*$ will be applied to (an inital
segment of) the last model of the phalanx ${\vec {\cal P}}$, ${\cal M}_\epsilon^{\cal
U}$. We may therefore finally
argue as in the proof of \cite[Theorem 8.6]{CMIP} to derive a contradiction.  
$\square$ {\scriptsize 
Claim 2}

\bigskip
Now let $(\kappa_\gamma \colon
\gamma < \theta)$ enumerate the cardinals of ${\bar K}$ in the
half-open interval
$[\eta,{\bar \lambda})$, and let $\lambda_\gamma = \kappa_\gamma^{+{\bar K}}$ for
$\gamma < \theta$. For each $\gamma < \theta$
we let $\delta(\gamma) < lh({\cal T})$ be the least $\delta$ such that 
${\cal M}^{\cal T}_{\delta} | \lambda_\gamma = {\bar K} | \lambda_\gamma$, 
we let ${\cal P}_\gamma$ be the largest initial segment ${\cal P}$ of 
${\cal M}^{\cal T}_{\delta(\gamma)}$ 
such that all bounded subsets of $\lambda_\gamma$ which are in ${\cal P}$ are in 
${\bar K}$ as well, and we
let $n(\gamma) < \omega$ be the $n<\omega$ such that
$$\rho_n({\cal P}_i) \leq \kappa_\gamma < \rho_{n-1}({\cal P}_i).$$
Notice that $n(\gamma)$ will always be defined by Claims 1 and 2. 
The following is easy to verify.

\bigskip
{\bf Claim 3.} $\{ n(\gamma) \ \colon \ \gamma < \theta \}$ is finite.

\bigskip
We now let for each $\gamma < \theta$,
$${\cal Q}_\gamma = Ult_{n-1}({\cal P}_\gamma;\pi \upharpoonright
{\bar K} | \lambda_\gamma).$$ 
Let $\sigma_\gamma$ denote the canonical embedding from ${\cal P}_\gamma$ into
${\cal Q}_\gamma$.
Unfortunately, again even if it is well-founded,
${\cal Q}_\gamma$ might not be premouse but rather a proto-mouse; 
this will in fact be
the case if ${\cal P}_\gamma$ has a top extender with critical point $\mu \leq
\kappa_\gamma$ and $n=1$, 
as then $\sigma_\gamma$ is discontinuous at $\mu^{+{\bar K}} =
\mu^{+{\cal P}_\gamma}$.

Let us therefore define, inductively 
for $\gamma < \theta$, objects ${\cal R}_\gamma$ as follows. We
understand that we let the construction break down as soon as one of the models
defined is ill-founded. Simultaneously to their definition we notice that
$\rho_\omega({\cal R}_\gamma) \leq \kappa_\gamma$, and we define $n^*(\gamma)$ as the
least $n<\omega$ such that $\rho_n({\cal R}_\gamma) \leq \kappa_\gamma
< \rho_{n-1}({\cal R}_\gamma)$.
We set ${\cal R}_\gamma = {\cal Q}_\gamma$ if
${\cal Q}_\gamma$ is a premouse. Otherwise, if $\kappa_{\bar \gamma}$ is the critical
point of the top extender of ${\cal P}_\gamma$,
we let $${\cal R}_\gamma = Ult_{n^*({\bar \gamma})}({\cal R}_{\bar \gamma};G){\rm , }$$
where $G$ is the top extender fragment of ${\cal Q}_\gamma$.

It is straightforward to see that if all the ${\cal R}_\gamma$ are well-defined 
then $\{ n^*(\gamma) \colon \gamma < \theta \} \subset 
\{ n(\gamma) \colon \gamma < \theta \}$, so that by Claim 3 we immediately get:

\bigskip
{\bf Claim 4.} $\{ n^*(\gamma) \ \colon \ \gamma < \theta \}$ is finite.

\bigskip
We now aim to prove that all the ${\cal R}_\gamma$ are well-defined. In fact: 

\bigskip
{\bf Claim 5.} For each $\gamma < \theta$ do we have that ${\cal R}_\gamma
\triangleleft K$.

\bigskip
{\sc Proof.} Fix $\gamma < \theta$.
We assume that ${\cal Q}_\gamma$ is a premouse,
leaving the rest as an exercise for the reader.
Let ${\tilde \lambda} = sup \pi{\rm " } \lambda_\kappa$.
By standard arguments in oder to show Claim 5 it will certainly be enough to prove the
following.

\bigskip
{\bf Subclaim 2.} $((K,{\cal Q}_\gamma),{\tilde \lambda})$ is iterable.

\bigskip
{\sc Proof.} Let $\sigma \colon {\bar
H} \rightarrow H_\Omega$ be such that $\Omega$ is regular and large enough,
${\bar H}$ is countable and transitive, and $\{ ((K,{\cal Q}_\gamma),{\tilde \lambda}),
{\cal P}_\gamma \} 
\subset ran(\sigma)$. It suffices to prove that
$\sigma^{-1}(((K,{\cal Q}_\gamma),{\tilde \lambda}))$ can be embedded into
an iterable phalanx.

Let $n = n(\gamma)$. Let $${\cal P}^* = Ult_{n-1}(\sigma^{-1}({\cal P}_\gamma);
\sigma \upharpoonright \sigma^{-1}({\cal P}_\gamma | \lambda_\gamma)){\rm , }$$
and let $k$ be the canonical embedding from ${\cal P}^*$ into ${\cal P}_\gamma$.
As $cf(\lambda_\gamma) > \omega$, the critical point of $k$ will be
$k^{-1}(\lambda_\gamma)$. Using \cite[\S 8, Lemma 4]{ronald} we then get as above that
in fact ${\cal P}^* \triangleleft {\cal P}_\gamma | \lambda_\gamma = {\bar K} |
\lambda_\gamma$. We claim that
$$[a,f] \mapsto \pi \circ k^{-1}(f)(a)$$
defines a sufficiently elementary 
embedding from $ran(\sigma) \cap {\cal Q}_\gamma$ into $\pi({\cal P}^*)$.
Here, $a \in [{\tilde \lambda}]^{<\omega} \cap ran(\sigma)$ and
$f$ is a $\Sigma_1^{(n-1)}({\cal P}_\gamma)$ good function with parameters from
$ran(\sigma) \cap {\cal P}_\gamma$.
The reason for this is that for appropriate $a$ and $f$ do we have that
\begin{eqnarray*}
& ran(\sigma) \cap {\cal Q}_\gamma \models \Phi([a,f]) {\rm \ \ iff } & \\
& a \in \pi(\{ u \colon {\cal P}_\gamma \models \Phi(f(u)) \}) {\rm \ \ iff } & \\
& a \in \pi(\{ u \colon {\cal P}^* \models \Phi(k^{-1}(f)(u)) \}) {\rm \ \ iff } & \\
& a \in \{ u \colon \pi({\cal P}^*) \models \Phi(k^{-1}(f)(u)) \} {\rm \ \ iff } & \\
& \pi({\cal P}^*) \models \Phi(\pi \circ k^{-1}(f)(a). & 
\end{eqnarray*}
This now implies that we my embedded
$\sigma^{-1}(((K,{\cal Q}_\gamma),{\tilde \lambda}))$ into
the phalanx $$((K,\pi({\cal P}^*)) , k^{-1}(\lambda_\gamma)).$$ 
However, $\pi({\cal P}^*)
\triangleleft K$, and hence this latter phalanx is iterable.
$\square$ {\scriptsize Subclaim 2} $\square$ {\scriptsize Claim 5} 

\bigskip
We have shown that if $i \in N \cap \lambda$, say
$\sigma(\kappa_\gamma) = \aleph_{i+1}$, then $n(sup N \cap \aleph_{i+1}) = 
n^*(\gamma)$.
By Claim 4 this establishes the Theorem.
$\square$ {\scriptsize Theorem \ref{main-thm}} 

\bigskip
Many questions remain open. Can Theorem \ref{main-thm} be extended to the core model
of \cite{CMIP}? Or can Theorem \ref{main-thm} even be extended to fine structural
models which do not know how to fully iterate themselves? 
Finally: Can one use methods provided by the current
paper to get a reasonable lower bound for the consistency strength of the assumption
that $(S_n \colon n<\omega)$ must be mutually stationary provided 
every individual $S_n \subset
\aleph_n$ is stationary?

\end{document}